\documentclass[reqno]{amsart}
\usepackage{amssymb,hyperref}

\theoremstyle{plain}
\newtheorem{thm}{Theorem}
\newtheorem{lem}{Lemma}

\theoremstyle{definition}
\newtheorem{ex}{Example}

\renewcommand{\Re}{\mathrm{Re}}

\title{An application of Miller and Mocanu lemma}

\author{Hitoshi Shiraishi}
\address{Hitoshi Shiraishi \newline
Department of Mathematics \newline
Kinki University \newline
Higashi-Osaka, Osaka 577-8502, Japan}
\email{shiraishi@math.kindai.ac.jp}

\author{Shigeyoshi Owa}
\address{Shigeyoshi Owa \newline
Department of Mathematics \newline
Kinki University \newline
Higashi-Osaka, Osaka 577-8502, Japan}
\email{owa@math.kindai.ac.jp}


\date{}

\begin{document}

\begin{abstract}
Let $\mathcal{H}[a_0,n]$ be the class of functions $f(z)=a_0+a_nz^n+\ldots$
which are analytic in the open unit disk $\mathbb{U}$.
For $f(z)\in\mathcal{H}[a_0,n]$,
S. S. Miller and P. T. Mocanu (J. Math. Anal. Appl. {\bf 65}(1978), 289-305)
have shown Miller and Mocanu lemma
which is the generalization of Jack lemma by I. S. Jack (J. London Math. Soc. {\bf 3}(1971), 469-474).
Applying Miller and Mocanu lemma,
an interesting property for $f(z)\in\mathcal{H}[a_0,n]$ and an example are discuss.
\end{abstract}

\begin{flushleft}
This paper was published in the journal: \\
Stud. Univ. Babe\c{s}-Bolyai Math. {\bf 55} (2010), No. 3, 207--211. \\
\url{http://www.cs.ubbcluj.ro/~studia-m/2010-3/shiraishi-final.pdf}
\end{flushleft}
\hrule

\

\

\maketitle

\section{Introduction}

\

Let $\mathcal{H}[a_0,n]$ denote the class of functions $f(z)$ of the form
$$
f(z)
= a_0 + \sum_{k=n}^{\infty} a_k z^k
\qquad (n=1,2,3,\ldots)
$$
which are analytic in the open unit disk $\mathbb{U}=\{ z\in\mathbb{C}:|z|<1 \}$,
where $a_0\in\mathbb{C}$.
Jack \cite{m1ref1} has shown the result for analytic functions $w(z)$ in $\mathbb{U}$ with $w(0)=0$, which is called Jack's lemma.
In 1978,
Miller and Mocanu \cite{m1ref2} have given the generalization theorem for Jack's lemma, 
which was called Miller and Mocanu lemma.

\begin{lem}[Miller and Mocanu lemma] \label{d5lem1} \quad
Let $f(z) \in \mathcal{H}[a_0,n]$
with $f(z) \not\equiv a_0$.
If there exists a point $z_0 \in \mathbb{U}$ such that
$$
\max_{|z| \leqq |z_0|} |f(z)|
=|f(z_0)|,
$$
then
$$
\frac{z_0 f'(z_0)}{f(z_0)}
= m
$$
and
$$
\Re \left( \frac{z_0 f''(z_0)}{f'(z_0)} \right)+1
\geqq m,
$$
where $m$ is real and
$$
m
\geqq n \frac{|f(z_0)-a_0|^2}{|f(z_0)|^2-|a_0|^2}
\geqq n \frac{|f(z_0)|-|a_0|}{|f(z_0)|+|a_0|}.
$$
\end{lem}

If $a_0=0$,
then the above lemma becomes Jack's lemma due to Jack \cite{m1ref1}.

\

\section{Main theorem}

\

Applying Miller and Mocanu lemma,
we derive

\begin{thm} \label{d5thm1} \quad
Let $f(z) \in \mathcal{H}[a_0,n]$ with $f(z) \neq 0$ for $z\in\mathbb{U}$.
If there exists a point $z_0\in\mathbb{U}$ such that
$$
\min_{|z| \leqq |z_0|} |f(z)|
=|f(z_0)|,
$$
then
\begin{equation}
\frac{z_0 f'(z_0)}{f(z_0)}
= -m
\label{d5thm1eq1}
\end{equation}
and
\begin{equation}
\Re\left(\frac{z_0 f''(z_0)}{f'(z_0)}\right)+1
\geqq -m,
\label{d5thm1eq2}
\end{equation}
where
$$
m
\geqq n \frac{|a_0-f(z_0)|^2}{|a_0|^2-|f(z_0)|^2}
\geqq n \frac{|a_0|-|f(z_0)|}{|a_0|+|f(z_0)|}.
$$
\end{thm}

\

\begin{proof} \quad
We defined the function $g(z)$ by
\begin{align*}
g(z)
& = \frac{1}{f(z)} \\
& = c_0 + c_n z^n + c_{n+1} z^{n+1} + \ldots
\qquad \left( c_0=\frac{1}{a_0} \right).
\end{align*}

Then,
$g(z)$ is analytic in $\mathbb{U}$ and $g(0)=c_0 \ne 0$.
Furthermore,
by the assumtion of the theorem,
$|g(z)|$ takes its maximum value at $z=z_0$
in the closed disk $|z| \leqq |z_0|$.
It follows from this that
$$
|g(z_0)|
= \frac{1}{|f(z_0)|}
= \frac{1}{\displaystyle{\min_{|z| \leqq |z_0|}} |f(z)|}
= \max_{|z| \leqq |z_0|} |g(z)|.
$$

Therefore,
applying Lemma \ref{d5lem1} to $g(z)$,
we observe that
$$
\frac{z_0 g'(z_0)}{g(z_0)}
= - \frac{z_0 f'(z_0)}{f(z_0)}
= m
$$
which shows (\ref{d5thm1eq1}) and
\begin{align*}
\Re\left(\frac{z_0 g''(z_0)}{g'(z_0)}\right) +1
& = \Re \left( \frac{z_0f''(z_0)}{f'(z_0)} -2\frac{z_0f'(z_0)}{f(z_0)} \right) +1 \\
&= \Re\left(\frac{z_0 f''(z_0)}{f'(z_0)}\right) +2m +1 \\
&\geqq m
\end{align*}
which implies (\ref{d5thm1eq2}),
where
$$
m
\geqq n \frac{|g(z_0)-c_0|^2}{|g(z_0)|^2-|c_0|^2}
= n \frac{|a_0-f(z_0)|^2}{|a_0|^2-|f(z_0)|^2}
\geqq n \frac{|a_0|-|f(z_0)|}{|a_0|+|f(z_0)|}.
$$

This completes the asserion of Theorem \ref{d5thm1}.
\end{proof}

\

\begin{ex} \label{d5ex1} \quad
Let us consider the function $f(z)$ given by
\begin{align*}
f(z)
&= \frac{a_0+ \left( e^{i\arg(a_0)}-a_0 \right) z^n}{1-z^n} \\
&= a+e^{i\arg(a_0)}z^n+e^{i\arg(a_0)}z^{2n}+\ldots
\qquad (z\in\mathbb{U})
\end{align*}
for some complex number $a_0$ with $|a_0| > \dfrac{1}{2}$.
Then,
$f(z)$ maps the disk $\mathbb{U}_r= \{ z:|z|<r\leqq1 \}$
onto the domain
$$
\left| f(z)- \left( a_0+\frac{e^{i\arg(a_0)}r^{2n}}{1-r^{2n}} \right) \right|
\leqq \frac{r^n}{1-r^{2n}}.
$$

Thus,
we know that there exists a point $z_0=re^{i\frac{\pi}{n}}\in\mathbb{U}$ such that
$$
\min_{|z|\leqq|z_0|}|f(z)|
= |f(z_0)|
= |a_0|-\frac{r^n}{1+r^n}.
$$

For such a point $z_0$,
we obtain that
$$
\frac{z_0 f'(z_0)}{f(z_0)}
= -\frac{nr^n}{(1+r^n)(|a_0|-(1-|a_0|)r^n)}
= -m
$$
where
$$
m
= \frac{nr^n}{(1+r^n)(|a_0|-(1-|a_0|)r^n)}
> 0.
$$

Therefore,
we get that
$$
\Re\left(\frac{z_0 f''(z_0)}{f'(z_0)}\right)+1
= n\frac{1-r^n}{1+r^n}
> 0
> -m.
$$

Furthermore,
we obtain that
$$
n\frac{|a_0-f(z_0)|^2}{|a_0|^2-|f(z_0)|^2}
= \frac{nr^n}{2|a_0|+(2|a_0|-1)r^n}
= \frac{nr^n}{2 \left( |a_0|-(1-|a_0|)r^n+\dfrac{1}{2}r^n \right)}
< m.
$$
\end{ex}

\

Putting $a_0$ with a real number in Example \ref{d5ex1},
we get Example \ref{d5ex2}.

\

\begin{ex} \label{d5ex2} \quad
Let us consider the function
\begin{align*}
f(z)
&= \frac{a_0+(1-a)z^n}{1-z^n} \\
&= a_0+z^n+z^{2n}+\ldots
\qquad (z\in\mathbb{U})
\end{align*}
for $a_0 > \dfrac{1}{2}$.
Then,
it follows that the function $f(z)$ maps the disk $\mathbb{U}_r$
onto the domain
$$
\left| f(z)- \left( a_0+\frac{r^{2n}}{1-r^{2n}} \right) \right|
\leqq \frac{r^n}{1-r^{2n}}.
$$

Thus,
there exists a point $z_0=re^{i\frac{\pi}{n}}\in\mathbb{U}$ such that
$$
\min_{|z|\leqq|z_0|}|f(z)|
= |f(z_0)|
= a_0-\frac{r^n}{1+r^n}.
$$

For such a point $z_0$,
we obtain
$$
\frac{z_0 f'(z_0)}{f(z_0)}
= -\frac{nr^n}{(1+r^n)(a_0-(1-a_0)r^n)}
= -m
$$
where
$$
m
= \frac{nr^n}{(1+r^n)(a_0-(1-a_0)r^n)}
> 0.
$$

Therefore,
we see that
$$
\Re\left(\frac{z_0 f''(z_0)}{f'(z_0)}\right)+1
= n\frac{1-r^n}{1+r^n}
> 0
> -m.
$$

Moreover,
we have that
$$
n\frac{|a_0-f(z_0)|^2}{|a_0|^2-|f(z_0)|^2}
= \frac{nr^n}{2a_0+(2a_0-1)r^n}
= \frac{nr^n}{2 \left( a_0-(1-a_0)r^n+\dfrac{1}{2}r^n \right)}
< m.
$$
\end{ex}

\

\end{document}